# On a family that unifies Generalized Marshall-Olkin and Poisson-G family of distribution

[1*] Laba Handique, [2] Farrukh Jamal and [3] Subrata Chakraborty

[1, 3] Department of Statistics, Dibrugarh University, Dibrugarh-786004, India
[2] Department of Statistics, Government S.A. Post-Graduate College, Dera Nawab Sahib, Pakistan.

**Abstract**

Unifying the generalized Marshall-Olkin (GMO) and Poisson-G (P-G) a new family of distribution is proposed. Density and the survival function are expressed as infinite mixtures of P-G family. The quantile function, asymptotes, shapes, stochastic ordering, moment generating function, order statistics, probability weighted moments and Rényi entropy are derived. Maximum likelihood estimation with large sample properties is presented. A Monte Carlo simulation is used to examine the pattern of the bias and the mean square error of the maximum likelihood estimators. An illustration of comparison with some of the important sub models of the family in modeling a real data reveals the utility of the proposed family.

______________________________________________

## 1. Introduction

With the basic motivation to bring in more flexibility in the modelling different type of data, a preferred area of research in the field of probability distribution is that of generating new distributions starting with a base line distribution by inducing one or more additional parameters through various methodologies. A number of useful continuous univariate-G families have been added in the literature in recent times. Notable families introduced since 2017 are Poisson-G family (Abouelmagd *et al*., 2017), Beta-G Poisson family (Gokarna *et al*., 2018), Marshall-Olkin Kumaraswamy-G family (Handique *et al*., 2017), Generalized Marshall-Olkin Kumaraswamy-G family (Chakraborty and Handique 2017), Exponentiated generalized-G Poisson family (Gokarna and Haitham, 2017), Beta Kumaraswamy-G family (Handique *et al*., 2017), Beta generated Kumaraswamy Marshall-Olkin-G family (Handique and Chakraborty, 2017a), Beta generalized Marshall-Olkin Kumaraswamy-G family (Handique and Chakraborty, 2017b), Beta generated Marshall-Olkin-Kumaraswamy-G (Chakraborty *et al*., 2018), Exponentiated generalized Marshall-Olkin-G family by (Handique *et al*., 2019), Kumaraswamy generalized Marshall-Olkin-G family (Chakraborty and Handique, 2018), Odd modified exponential generalized family (Ahsan *et al*., 2018), Zografos-Balakrishnan Burr XII family (Emrah *et al.,* 2018), Zero truncated Poisson family (Abouelmagd *et al*., 2019), Extended generalized Gompertz family (Thiago *et al*., 2019), Odd Half-Cauchy family (Chakraborty *et al*., 2020), Generalized modified exponential-G family (Handique *et al*., 2020), Beta Poisson-G family (Handique *et al*., 2020) and Kumaraswamy Poisson-G (Chakraborty *et al*., 2020) among others.

In this current article, we propose a new family of continuous probability distribution to unify generalized Marshall-Olkin (GMO) of Jayakumar and Mathew (2008) and the Poisson-G (P-G) family of distribution (Tahir *et al*., 2016) and call it the Generalized Marshall-Olkin Poisson-G ($\text{GMOP} - \text{G}(\theta, \alpha, \lambda)$). Now we briefly describe the GMO and P-G family and then introduce GMOP-G in the next section.

**GMO**: Jayakumar and Mathew (2008) proposed a generalization of the Marshall-Olkin family (Marshall and Olkin, 1997) called generalized Marshall-Olkin (GMO) family of distribution with baseline distribution having survival function (sf) $\overline{F}(t)$ and probability distribution function (pdf) $f(t)$. The sf and pdf of the GMO distribution are given respectively by

$$\overline{F}^{GMO}(t;\theta,\alpha) = \left[\frac{\alpha \overline{F}(t)}{1-\overline{\alpha} \overline{F}(t)}\right]^{\theta} \text{ and } f^{GMO}(t;\theta,\alpha) = \frac{\theta \alpha^{\theta} f(t) \overline{F}(t)^{\theta-1}}{[1-\overline{\alpha} \overline{F}(t)]^{\theta+1}}, \tag{1}$$

where $-\infty < t < \infty$, $\alpha > 0$ ($\overline{\alpha} = 1-\alpha$) and $\theta > 0$ is an additional shape parameter. When $\theta = 1$, $\overline{F}^{GMO}(t;\theta,\alpha) = \overline{F}^{MO}(t;\alpha)$ and for $\alpha = \theta = 1$, $\overline{F}^{GMO}(t;\theta,\alpha) = \overline{F}(t)$.

**P-G**: The Poisson-G (P-G) family of distribution with sf and cdf is given by (Chakraborty *et al.*, 2020)

$$\overline{G}^{PG}(t;\lambda) = \frac{e^{-\lambda G(t)} - e^{-\lambda}}{1-e^{-\lambda}} \text{ and } G^{PG}(t;\lambda) = \frac{1-e^{-\lambda G(t)}}{1-e^{-\lambda}}, \lambda \in R - \{0\}; n = 1, 2,... \tag{2}$$

where $G(t)$ is the baseline distribution. The corresponding pdf of (2) is given by

$$g^{PG}(t;\lambda) = (1-e^{-\lambda})^{-1} \lambda g(t) e^{-\lambda G(t)}, \lambda \in R - \{0\}; -\infty < t < \infty. \tag{3}$$

Rest of the article is arranged in 5 more sections. In Section 2 we introduce the proposed family along with its physical basis and list of some important sub models also defined some mathematical properties. In Section 3, linear representation of the sf and pdf of the proposed family also we discuss some statistical properties of the proposed family. In Section 4, maximum likelihood methods of estimation of parameters and simulation are presented. The data fitting applications is presented in Section 5. The article ends with a conclusion in Section 6.

## 2. Generalized Marshall-Olkin Poisson-G

In this section we introduce the $\text{GMOP} - G(\theta, \alpha, \lambda)$ family and also provide its special cases and a statistical genesis.

The sf, cdf, pdf and hrf of this proposed distribution are respectively given by:

$$\overline{F}^{GMOPG}(t;\theta,\alpha,\lambda) = \left[\frac{\alpha(e^{-\lambda G(t)} - e^{-\lambda})}{1-\alpha e^{-\lambda} - \overline{\alpha} e^{-\lambda G(t)}}\right]^{\theta}, F^{GMOPG}(t;\theta,\alpha,\lambda) = 1 - \left[\frac{\alpha(e^{-\lambda G(t)} - e^{-\lambda})}{1-\alpha e^{-\lambda} - \overline{\alpha} e^{-\lambda G(t)}}\right]^{\theta}, \tag{4}$$

$$f^{GMOPG}(t;\theta,\alpha,\lambda) = \frac{\theta \lambda \alpha^{\theta} (1-e^{-\lambda}) g(t) e^{-\lambda G(t)} (e^{-\lambda G(t)} - e^{-\lambda})^{\theta-1}}{(1-\alpha e^{-\lambda} - \overline{\alpha} e^{-\lambda G(t)})^{\theta+1}}, \tag{5}$$

and

$$h^{GMOPG}(t;\theta,\alpha,\lambda) = \frac{\theta \lambda (1-e^{-\lambda}) g(t) e^{-\lambda G(t)} (e^{-\lambda G(t)} - e^{-\lambda})^{-1}}{1-\alpha e^{-\lambda} - \overline{\alpha} e^{-\lambda G(t)}}$$

In particular, we get for

(i) $\theta = 1$, the $\text{MOP-G}(\alpha, \lambda)$ distribution.

(ii) $\theta = \alpha = 1$, the $P - G(\lambda)$ distribution.

(iii) $\lambda = 0$, the $\text{GMO}(\theta, \alpha)$ distribution.

(iv) $\theta = 1$, $\lambda = 0$, the $\text{MO}(\alpha)$ distribution.

**Proposition 1:** Let $T_{i1}, T_{i2}, ..., T_{iN}$, $i = 1, 2, \cdots, \theta$ be a sequence of $\theta N$ i.i.d. random variables from Poisson-G distribution and $W_i = \min(T_{i1}, T_{i2}, ..., T_{iN})$ and $V_i = \max(T_{i1}, T_{i2}, ..., T_{iN})$. Then

(i) $\min_i W_i$ follows $\text{GMOP-G}(\theta, \alpha, \lambda)$ if $Geometric(\alpha)$

(ii) $\max_i V_i$ follows $\text{GMOP-G}(\theta, \alpha, \lambda)$ if $Geometric(1/\alpha)$.

**Proof:**

**Case** (i) When $0 < \alpha \leq 1$, considering $N$ has a geometric distribution with parameter $\alpha$, we get

$$P[\min\{W_1, W_2, ..., W_\theta\} > t] = P[W_1 > t] P[W_2 > t] ... P[W_\theta > t]$$

$$= \prod_{i=1}^{\theta} P[W_i > t] = [\bar{F}^{\text{MOPG}}(t; \alpha, \lambda)]^\theta = \left[ \frac{\alpha(e^{-\lambda G(t)} - e^{-\lambda})}{1 - \alpha e^{-\lambda} - \bar{\alpha} e^{-\lambda G(t)}} \right]^\theta.$$

**Case** (ii) For $\alpha > 1$, considering $N$ has a geometric distribution with parameter $1/\alpha$, we get

$$P[\min\{V_1, V_2, ..., V_\theta\} > t] = P[V_1 > t] P[V_2 > t] ... P[V_\theta > t]$$

$$= \prod_{i=1}^{\theta} P[V_i > t] = [\bar{F}^{\text{MOPG}}(t; \alpha, \lambda)]^\theta = \left[ \frac{\alpha(e^{-\lambda G(t)} - e^{-\lambda})}{1 - \alpha e^{-\lambda} - \bar{\alpha} e^{-\lambda G(t)}} \right]^\theta \quad \square$$

In what follows we investigate some general properties, parameter estimation and real life applications.

## 2.1 Special models and shape of the density and hazard function

In this section we have plotted the pdf and hrf of the $\text{GMOP-G}(\theta, \alpha, \lambda)$ taking $G$ to be Weibull (W) and exponential (E) for some chosen values of the parameters to show the variety of shapes assumed by the family.

The pdf and hrf of these distributions are obtained from $\text{GMOP-G}(\theta, \alpha, \lambda)$ as follows:

➢ The GMOP-Weibull (GMOP-W) distribution

Considering the Weibull distribution (Weibull, 1951) with parameters $\beta > 0$ and $\delta > 0$ having pdf and cdf $g(t) = \delta \beta t^{\delta - 1} e^{-\beta t^\delta}$ and $G(t) = 1 - e^{-\beta t^\delta}$ respectively we get the pdf and hrf of $\text{GMOP-W}(\theta, \alpha, \lambda, \beta, \delta)$ distribution as

$$f^{\text{GMOPW}}(t; \theta, \alpha, \lambda, \beta, \delta) = \frac{\theta \lambda \alpha^\theta (1 - e^{-\lambda}) \delta \beta t^{\delta - 1} e^{-\beta t^\delta} e^{-\lambda(1 - e^{-\beta t^\delta})} (e^{-\lambda(1 - e^{-\beta t^\delta})} - e^{-\lambda})^{\theta - 1}}{(1 - \alpha e^{-\lambda} - \bar{\alpha} e^{-\lambda(1 - e^{-\beta t^\delta})})^{\theta + 1}},$$

and $\quad h^{\text{GMOPW}}(t; \theta, \alpha, \lambda, \beta, \delta) = \dfrac{\theta \lambda \alpha^\theta (1 - e^{-\lambda}) \delta \beta t^{\delta - 1} e^{-\beta t^\delta} e^{-\lambda(1 - e^{-\beta t^\delta})} (e^{-\lambda(1 - e^{-\beta t^\delta})} - e^{-\lambda})^{-1}}{1 - \alpha e^{-\lambda} - \bar{\alpha} e^{-\lambda(1 - e^{-\beta t^\delta})}}.$

Taking $\delta = 1$ in $\text{GMOP-W}(\theta, \alpha, \lambda, \beta, \delta)$ we get the $\text{GMOP-E}(\theta, \alpha, \lambda, \beta)$ with pdf and hrf is given by

$$f^{\text{GMOPE}}(t; \theta, \alpha, \lambda, \beta) = \frac{\theta \lambda \alpha^\theta (1 - e^{-\lambda}) \beta e^{-\beta t} e^{-\lambda(1 - e^{-\beta t})} (e^{-\lambda(1 - e^{-\beta t})} - e^{-\lambda})^{\theta - 1}}{(1 - \alpha e^{-\lambda} - \bar{\alpha} e^{-\lambda(1 - e^{-\beta t})})^{\theta + 1}},$$

and $$h^{\text{GMOPE}}(t;\theta,\alpha,\lambda,\beta) = \frac{\theta\lambda\alpha^{\theta}(1-e^{-\lambda})\beta e^{-\beta t}e^{-\lambda(1-e^{-\beta t})}(e^{-\lambda(1-e^{-\beta t})}-e^{-\lambda})^{-1}}{1-\alpha e^{-\lambda}-\overline{\alpha}e^{-\lambda(1-e^{-\beta t})}}.$$

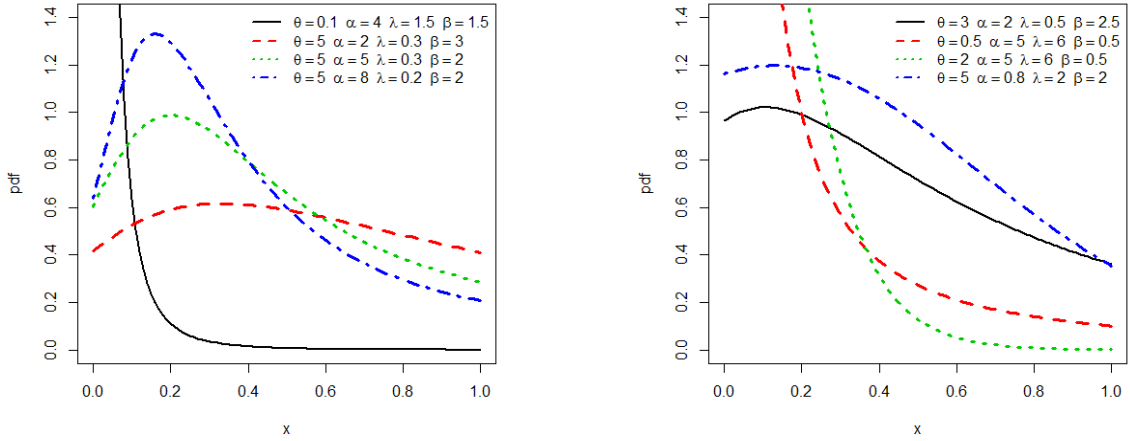

**Fig 1: pdf plots of the $\text{GMOP-E}(\theta,\alpha,\lambda,\beta)$ distribution**

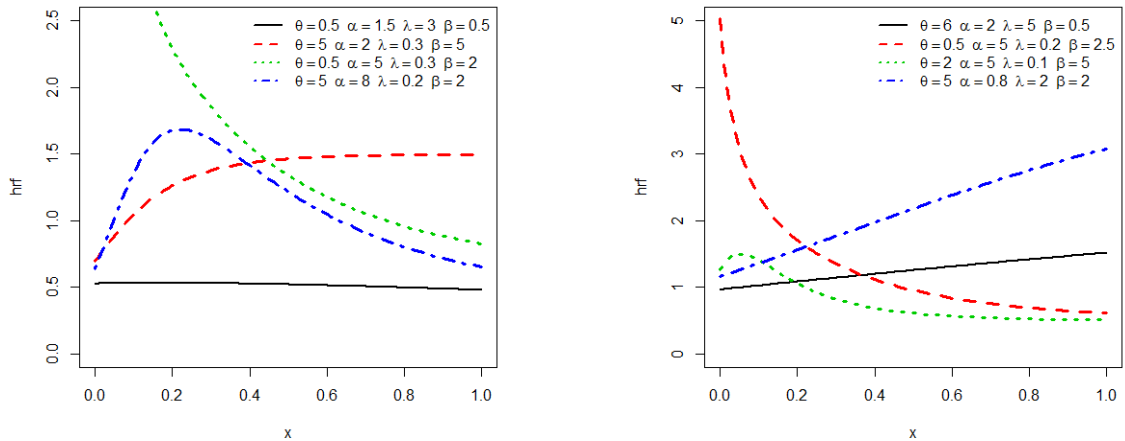

**Fig 2: hrf plots of the $\text{GMOP-E}(\theta,\alpha,\lambda,\beta)$ distribution**

### 2.2 Quantile and related measures

The $p^{th}$ quantile $t_p$ for $\text{GMOP}-\text{G}(\theta,\alpha,\lambda)$ can be easily obtained by solving the equation $F^{\text{GMOPG}}(t)=p$ as

$$t_p = G^{-1}\left[-\frac{1}{\lambda}\log\left[\frac{\alpha e^{-\lambda}+(1-\alpha e^{-\lambda})(1-F(t))^{1/\theta}}{\alpha+\overline{\alpha}(1-F(t))^{1/\theta}}\right]\right].$$

A random number '$t$' from $\text{GMOP}-\text{G}(\theta,\alpha,\lambda)$ via an uniform random number '$u$' can be generated by using the formula

$$t = G^{-1}\left[-\frac{1}{\lambda}\log\left[\frac{\alpha e^{-\lambda} + (1-\alpha e^{-\lambda})(1-u)^{1/\theta}}{\alpha + \bar{\alpha}(1-u)^{1/\theta}}\right]\right]. \qquad (6)$$

For example, when we consider the exponential distribution having pdf and cdf as $g(x:\beta) = \beta e^{-\beta x}$, $x > 0$ $\beta > 0$, and $G(x:\beta) = 1 - e^{-\beta x}$, the $p^{th}$ quantile, $t_p$ is given by

$$t_p = -\frac{1}{\beta}\log\left[1 + \frac{1}{\lambda}\log\left\{\frac{\alpha e^{-\lambda} + (1-\alpha e^{-\lambda})(1-u)^{1/\theta}}{\alpha + \bar{\alpha}(1-u)^{1/\theta}}\right\}\right].$$

Here the flexibility of skewness and kurtosis of $\text{GMOP} - \text{G}(\theta, \alpha, \lambda)$ is checked by plotting Galton skewness (S) that measures the degree of the long tail and Morris (1988) kurtosis (K) that measures the degree of tail heaviness. These are respectively defined by

$$S = \frac{Q(6/8) - 2Q(4/8) + Q(2/8)}{Q(6/8) - Q(2/8)} \text{ and } K = \frac{Q(7/8) - Q(5/8) + Q(3/8) - Q(1/8)}{Q(6/8) - Q(2/8)}.$$

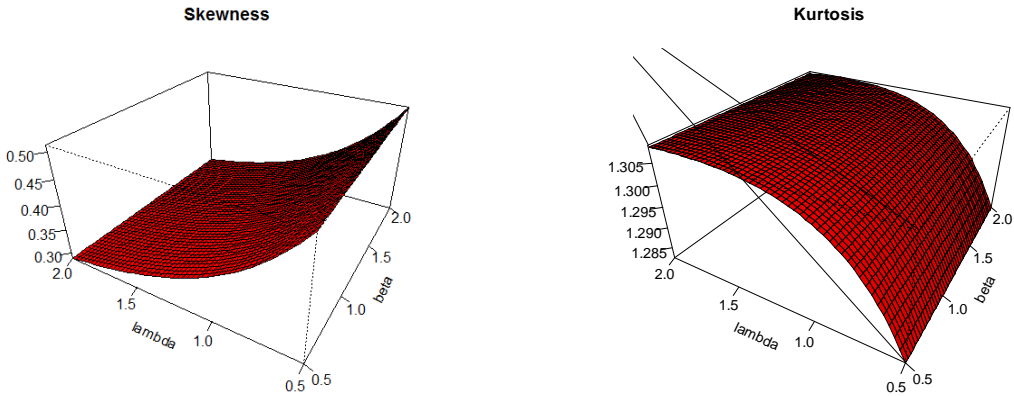

**Fig 3:** Plots of the Galton skewness $S$ and the Moor kurtosis $K$ for the GMOP-E distribution with parameters $\theta = 3, \alpha = 2, 0.2 < \lambda, \beta < 2$

## 2.3 Asymptotes and shapes

Two propositions regarding asymptotes of the proposed family are discussed here.

**Proposition 2:** The asymptotes of pdf, cdf and hrf of $\text{GMOP} - \text{G}(\theta, \alpha, \lambda)$ as $t \to 0$ are given by

$$f^{\text{GMOPG}}(t; \theta, \alpha, \lambda) \sim \theta \lambda g(t) / \alpha (1 - e^{-\lambda}),$$

$$F^{\text{GMOPG}}(t; \theta, \alpha, \lambda) \sim 0 \qquad \text{and}$$

$$h^{\text{GMOPG}}(t; \theta, \alpha, \lambda) \sim \theta \lambda g(t) / \alpha (1 - e^{-\lambda}).$$

**Proposition 3:** The asymptotes of pdf, cdf and hrf of $\text{GMOP} - \text{G}(\theta, \alpha, \lambda)$ as $t \to \infty$ are given by

$$f^{\text{GMOPG}}(t;\theta,\alpha,\lambda) \sim \theta\lambda\alpha^{\theta}e^{-\lambda}\,g(t)\,(e^{-\lambda G(t)}-e^{-\lambda})^{\theta-1}/(1-e^{-\lambda})^{\theta},$$

$$F^{\text{GMOPG}}(t;\theta,\alpha,\lambda) \sim 1-\alpha^{\theta}(e^{-\lambda G(t)}-e^{-\lambda})^{\theta}/(1-e^{-\lambda})^{\theta} \text{ and}$$

$$h^{\text{GMOPG}}(t;\theta,\alpha,\lambda) \sim \theta\lambda e^{-\lambda}\,g(t)\,(e^{-\lambda G(t)}-e^{-\lambda})^{-1}.$$

The shapes of the density and hazard rate function can be described analytically. The critical points of the $\text{GMOP}-\text{G}(\theta,\alpha,\lambda)$ family density function are the roots of the equation

$$\frac{g'(t)}{g(t)} - \lambda g(t) - (\theta-1)\frac{\lambda e^{-\lambda G(t)}g(t)}{e^{-\lambda G(t)}-e^{-\lambda}} - (\theta+1)\frac{\bar{\alpha}\lambda e^{-\lambda G(t)}g(t)}{1-\alpha e^{-\lambda}-\bar{\alpha}e^{-\lambda G(t)}} = 0. \quad (7)$$

The critical point of $\text{GMOP}-\text{G}(\theta,\alpha,\lambda)$ family hazard rate are the roots of the equation

$$\frac{g'(t)}{g(t)} - \lambda g(t) + \frac{\lambda e^{-\lambda G(t)}g(t)}{e^{-\lambda G(t)}-e^{-\lambda}} - \frac{\bar{\alpha}\lambda e^{-\lambda G(t)}g(t)}{1-\alpha e^{-\lambda}-\bar{\alpha}e^{-\lambda G(t)}} = 0. \quad (8)$$

There may be more than one root of (7) and (8). If $t=t_0$ is a root then it is a local maximum, a local minimum or a point of inflexion if $\psi(t_0)<0$, $\psi(t_0)>0$ or $\psi(t_0)=0$ and for (8) if $\omega(t_0)<0$, $\omega(t_0)>0$ or $\omega(t_0)=0$ where $\psi(t)=(d^2/dt^2)\log[f(t)]$ and $\omega(t)=(d^2/dt^2)\log[h(t)]$.

## 2.4 Stochastic orderings

Let $X$ and $Y$ be two random variables with cdfs $F$ and $G$, respectively, corresponding pdf's $f$ and $g$. Then $X$ is said to be smaller than $Y$ in the likelihood ratio order ($X \leq_{lr} Y$) if $f(t)/g(t)$ is decreasing in $t\geq 0$. Here we present a result of likelihood ratio ordering.

**Theorem 1:** Let $X \sim \text{GMOPG}(\theta,\alpha_1,\lambda)$ and $Y \sim \text{GMOPG}(\theta,\alpha_2,\lambda)$. If $\alpha_1<\alpha_2$, then $X \leq_{lr} Y$

Proof: $\dfrac{f(t)}{g(t)} = \left(\dfrac{\alpha_1}{\alpha_2}\right)^{\theta}\left[\dfrac{1-\alpha_2 e^{-\lambda}-\bar{\alpha}_2 e^{-\lambda G(t)}}{1-\alpha_1 e^{-\lambda}-\bar{\alpha}_1 e^{-\lambda G(t)}}\right]^{\theta+1}$

$$\frac{d}{dt}(f(t)/g(t)) = (\theta+1)(\alpha_1/\alpha_2)^{\theta}(\alpha_1-\alpha_2)\frac{[1-\alpha_2 e^{-\lambda}-\bar{\alpha}_2 e^{-\lambda G(t)}]^{\theta}\lambda e^{-\lambda G(t)}g(t)(1-e^{-\lambda})}{[1-\alpha_1 e^{-\lambda}-\bar{\alpha}_1 e^{-\lambda G(t)}]^{\theta+2}}.$$

Now this is always less than 0 since $\alpha_1<\alpha_2$. Hence, $f(t)/g(t)$ is decreasing in $t$. That is $X \leq_{lr} Y$

## 3 Linear representation

Linear representation of sf and pdf, etc in terms of corresponding functions of known distributions is an important tool for further mathematical properties. In this section we present some important results for proposed family.

### 3.1 Expansions of the survival and density functions as infinite linear mixture

Here we express the sf and pdf of the $\text{GMOP}-\text{G}(\theta,\alpha,\lambda)$ as infinite linear mixture of the corresponding functions of $\text{P}-\text{G}(\lambda)$ distribution.

Consider the series representation $(1-z)^{-k} = \sum_{j=0}^{\infty} \frac{\Gamma(k+j)}{\Gamma(k)j!} z^j = \sum_{j=0}^{\infty} \frac{(j+k-1)!}{(k-1)!j!} z^j$. (9)

This is valid for $|z| < 1$ and $k > 0$, where $\Gamma(.)$ is the gamma function.

Using equation (9) in equation (4), for $\alpha \in (0,1)$ we obtain

$$\overline{F}^{\text{GMOPG}}(t;\theta,\alpha,\lambda) = \alpha^\theta \{\overline{G}^{\text{PG}}(t;\lambda)\}^\theta \sum_{j=0}^{\infty} \frac{(j+\theta-1)!}{(\theta-1)!j!} (1-\alpha)^j \{\overline{G}^{\text{PG}}(t;\lambda)\}^j$$

$$= \sum_{j=0}^{\infty} \eta'_j [\overline{G}^{\text{PG}}(t;\lambda)]^{j+\theta} \qquad (10)$$

$$= \sum_{j=0}^{\infty} \eta'_j \overline{G}^{\text{PG}}(t;\lambda(j+\theta))$$

Differentiating in equation (10) with respect to '$t$' we get

$$f^{\text{GMOPG}}(t;\theta,\alpha,\lambda) = g^{\text{PG}}(t;\lambda) \sum_{j=0}^{\infty} \eta_j [\overline{G}^{\text{PG}}(t;\lambda)]^{j+\theta-1} \qquad (11)$$

$$= -\sum_{j=0}^{\infty} \eta'_j \frac{d}{dt} [\overline{G}^{\text{PG}}(t;\lambda(j+\theta))] \qquad (12)$$

$$= \sum_{j=0}^{\infty} \eta'_j g^{\text{PG}}(t;\lambda(j+\theta)),$$

where $\eta'_j = \eta'_j(\alpha) = \binom{j+\theta-1}{j}(1-\alpha)^j \alpha^\theta$, $\eta_j = \eta_j(\alpha) = (j+\theta)\eta'_j$.

### 3.2 Moment generating function

The moment generating function (mgf) of $\text{GMOP} - \text{G}(\theta,\alpha,\lambda)$ family can be easily expressed in terms of those of the exponentiated $\text{P-G}(\lambda)$ distribution using the results of Section 2.1. For example using equation (12) it can be seen that

$$M_T^{\text{GMOPG}}(s) = E[e^{sT}] = \int_{-\infty}^{\infty} e^{st} f^{\text{GMOPG}}(t) dt = -\int_{-\infty}^{\infty} e^{st} \sum_{j=0}^{\infty} \eta'_j \frac{d}{dt} [\overline{G}^{\text{PG}}(t;\lambda(j+\theta))] dt$$

$$= -\sum_{j=0}^{\infty} \eta'_j \int_{-\infty}^{\infty} e^{st} [\overline{G}^{\text{PG}}(t;\lambda(j+\theta))] dt = \sum_{j=0}^{\infty} \eta_j M_T^{\text{PG}}(s)$$

where $M_T^{\text{PG}}(s)$ is the mgf of a $\text{P-G}(\lambda)$ distribution.

**Table 1:** Mean, variance, skewness and kurtosis of the GMOP –E distribution with different values of $\theta, \alpha, \lambda$ and $\beta$

| $\theta$ | $\alpha$ | $\lambda$ | $\beta$ | Mean | Variance | Skewness | Kurtosis |
|---|---|---|---|---|---|---|---|
| 10 | 10 | 2 | 2 | 0.1572 | 0.0163 | 1.2096 | 4.6452 |
| 10 | 10 | 1 | 2 | 0.2224 | 0.0306 | 1.0761 | 4.0834 |
| 10 | 10 | 0.5 | 2 | 0.2648 | 0.0408 | 0.9749 | 3.7373 |
| 10 | 10 | 0.1 | 2 | 0.3033 | 0.0502 | 0.8840 | 3.4698 |
| 10 | 10 | 2 | 1 | 0.3145 | 0.0652 | 1.2096 | 4.6452 |
| 10 | 10 | 2 | 0.5 | 0.6291 | 0.2610 | 1.2096 | 4.6452 |
| 10 | 10 | 0.5 | 0.5 | 1.0592 | 0.6528 | 0.9749 | 3.7373 |
| 10 | 10 | 0.1 | 0.1 | 6.0675 | 20.1022 | 0.8840 | 3.4698 |
| 10 | 5 | 2 | 2 | 0.0918 | 0.0066 | 1.5144 | 6.0241 |
| 10 | 2 | 2 | 2 | 0.0419 | 0.0016 | 1.9171 | 8.4744 |
| 10 | 0.5 | 2 | 2 | 0.0115 | 0.0001 | 2.4163 | 12.8206 |
| 10 | 0.5 | 0.5 | 0.5 | 0.0834 | 0.0077 | 2.3537 | 12.1190 |
| 5 | 10 | 2 | 2 | 0.2692 | 0.0424 | 1.0978 | 4.3497 |
| 5 | 5 | 2 | 2 | 0.1676 | 0.0206 | 1.4448 | 5.8193 |
| 2 | 5 | 2 | 2 | 0.3615 | 0.0914 | 1.5105 | 6.3958 |
| 2 | 2 | 2 | 2 | 0.2049 | 0.0427 | 2.1968 | 10.7924 |
| 1 | 2 | 2 | 2 | 0.4180 | 0.19284 | 2.3136 | 11.2913 |
| 5 | 0.1 | 0.1 | 0.1 | 0.2309 | 0.0804 | 3.6320 | 30.5202 |
| 5 | 5 | 5 | 3 | 0.0489 | 0.0016 | 1.4027 | 5.7915 |
| 5 | 5 | 3 | 5 | 0.04882 | 0.0017 | 1.5071 | 6.2776 |
| 5 | 5 | 5 | 8 | 0.0183 | 0.0002 | 1.3967 | 5.7764 |
| 5 | 5 | 10 | 10 | 0.0069 | 0.00001 | 1.1168 | 2.1041 |

## 3.3 Rényi entropy

The entropy of a random variable is a measure of uncertainty variation and has been used in various situations in science and engineering. The Rényi entropy (see details, Song, 2001) is defined by

$$I_R(\delta) = (1-\delta)^{-1} \log\left(\int_{-\infty}^{\infty} f(t)^{\delta} dt\right) \quad , \text{ where } \delta > 0 \text{ and } \delta \neq 1.$$

Thus the Rényi entropy of $\text{GMOP} - \text{G}(\theta, \alpha, \lambda)$ distribution can be obtained as

$$I_R(\delta) = (1-\delta)^{-1} \log \left( \sum_{j=0}^{\infty} \mu_j \int_{-\infty}^{\infty} [g^{PG}(t;\lambda) \overline{G}^{PG}(t;\lambda)^{\theta-1}]^{\delta} [\overline{G}^{PG}(t;\lambda)]^j dt \right),$$

where $\mu_j = \mu_j(\alpha) = \{\theta^{\delta} \alpha^{\delta\theta} (1-\alpha)^j \Gamma[\delta(\theta+1)+j]\}/\{\Gamma[\delta(\theta+1)] \; j!\}$.

**Table 2:** Rényi entropy $GMOP-E(\theta, \alpha, \lambda, \beta)$ distribution with different values of $\theta, \alpha, \lambda$ and $\beta$

| Parameter | | | | $\delta$ | | | | | |
|---|---|---|---|---|---|---|---|---|---|
| $\theta$ | $\alpha$ | $\lambda$ | $\beta$ | **0.2** | **0.5** | **1.5** | **2** | **3** | **5** |
| 10 | 10 | 2 | 2 | -0.2550 | -0.6403 | -0.9916 | -1.0647 | -1.1549 | -1.2490 |
| 5 | 5 | 0.5 | 0.5 | 1.6816 | 1.3053 | 0.9722 | 0.9032 | 0.8176 | 0.7280 |
| 5 | 5 | 2 | 0.5 | 1.3469 | 0.8661 | 0.4519 | 0.3687 | 0.2669 | 0.1621 |
| 3 | 3 | 2 | 0.5 | 1.6090 | 1.0220 | 0.5252 | 0.42839 | 0.3113 | 0.1924 |
| 1.5 | 1.5 | 2 | 0.5 | 2.1288 | 1.3773 | 0.6945 | 0.5640 | 0.4097 | 0.2571 |
| 2 | 0.5 | 0.5 | 0.5 | 1.7182 | 0.8594 | 0.4390 | -0.1133 | -0.2982 | -0.4789 |

As expected the Rényi entropy turns out to be non increasing with $\delta$.

### 3.4 Order Statistics

Suppose $T_1, T_2, \ldots, T_n$ is a random sample from any $GMOP-G(\theta, \alpha, \lambda)$ distribution. Let $T_{i:n}$ denote the $i^{th}$ order statistics. The pdf of $T_{i:n}$ can be expressed as

$$f_{i:n}(t) = \frac{n!}{(i-1)!(n-i)!} f^{GMOPG}(t)[1-\overline{F}^{GMOPG}(t)]^{i-1} \overline{F}^{GMOPG}(t)^{n-i}$$

$$= \frac{n!}{(i-1)!(n-i)!} f^{GMOPG}(t) \sum_{l=0}^{i-1} (-1)^l \frac{(i-1)!}{l!(i-l-1)!} \overline{F}^{GMOPG}(t)^{n+l-i}$$

Now using the general expansion of the pdf and sf of $GMOP-G(\theta, \alpha, \lambda)$ distribution we get the pdf of the $i^{th}$ order statistics for of the $GMOP-G(\theta, \alpha, \lambda)$ as

$$f_{i:n}(t) = \frac{n!}{(i-1)!(n-i)!} [g^{PG}(t;\lambda) \sum_{j=0}^{\infty} \eta_j [\overline{G}^{PG}(t;\lambda)]^{j+\theta-1}]$$

$$\times \sum_{l=0}^{i-1} (-1)^l \binom{i-1}{l} [\sum_{k=0}^{\infty} \eta'_k [\overline{G}^{PG}(t;\lambda)]^{k+\theta(n+l-i)}], \text{ where } \eta_j \text{ and } \eta'_k \text{ are defined earlier.}$$

$$= g^{PG}(t;\lambda) \sum_{j,k=0}^{\infty} \Phi_{j,k} [\overline{G}^{PG}(t;\lambda)]^{j+k+\theta(n+l-i+1)-1} \tag{13}$$

$$= -\sum_{j,k=0}^{\infty} [\Phi_{j,k}/(j+k+\theta(n+l-i+1))] \frac{d}{dt} [\overline{G}^{PG}(t;\lambda)]^{j+k+\theta(n+l-i+1)}$$

$$= \sum_{j,k=0}^{\infty} \Phi'_{j,k} \, g^{PG}(t; \lambda(j+k+\theta(n+l-i+1)))],$$

where $\Phi_{j,k} = n\eta_j \eta'_k \binom{n-1}{i-1} \sum_{l=0}^{i-1} \binom{i-1}{l}(-1)^l$ and $\Phi'_{j,k} = \Phi_{j,k}/(j+k+\theta(n+l-i+1))$.

### 3.5 Probability weighted moments

The probability weighted moments (PWM), first proposed by Greenwood *et al.* (1979), are expectations of certain functions of a random variable whose mean exists. The $(p,q,r)^{th}$ PWM of $T$ is defined by

$$\Gamma_{p,q,r} = \int_{-\infty}^{\infty} t^p F(t)^q [1-F(t)]^r f(t)\,dt.$$

The $s^{th}$ moment of $T$ can be written as

$$E(T^s) = \sum_{j=0}^{\infty} \eta_j \int_0^{\infty} t^s [\overline{G}^{PG}(t;\lambda)]^{j+\theta-1} g^{PG}(t;\lambda)\,dt = \sum_{j=0}^{\infty} \eta_j \Gamma^{PG}_{s,0,j+\theta-1}$$ where $\eta_j$ define in Section 2.1 and

$$\Gamma_{p,q,r} = \int_0^{\infty} t^p \left(\frac{1-e^{-\lambda G(t)}}{1-e^{-\lambda}}\right)^q \left[1 - \frac{1-e^{-\lambda G(t)}}{1-e^{-\lambda}}\right]^r \frac{\lambda g(t) e^{-\lambda G(t)}}{1-e^{-\lambda}}\,dt \quad \text{is the PWM of } P\text{-}G(\lambda) \text{ distribution.}$$

Therefore the moments of the $\text{GMOP}-G(\theta,\alpha,\lambda)$ may be expressed in terms of the PWMs of $P\text{-}G(\lambda)$.

Similarly proceeding we can derive $s^{th}$ moment of the $i^{th}$ order statistic $T_{i:n}$, in a random sample of size $n$ from $\text{GMOP}-G(\theta,\alpha,\lambda)$ by using equations (13) as $E(T^s_{i,n}) = \sum_{j,k=0}^{\infty} \Phi_{j,k} \Gamma_{s,0,\,j+k+\theta(n+l-i+1)-1}$ where $\eta_j$ and $\Phi_{j,k}$ defined above.

### 4. Estimation

This section is devoted to the estimation of the $\text{GMOP}-G(\theta,\alpha,\lambda)$ model parameters via the maximum likelihood (ML) method.

### 4.1 Maximum likelihood method

Let $T = (t_1, t_2, \ldots, t_n)$ be a random sample of size $n$ from $\text{GMOP}-G(\theta,\alpha,\lambda)$ with parameter vector $\boldsymbol{\rho} = (\theta, \alpha, \lambda, \boldsymbol{\xi})$, where $\boldsymbol{\xi} = (\xi_1, \xi_2, \ldots, \xi_q)$ is the parameter vector of $G$. Then the log-likelihood function for $\boldsymbol{\rho}$ is given by

$$\ell = \ell(\boldsymbol{\rho}) = n\log(\theta\lambda\alpha^\theta) - n\log(1-e^{-\lambda}) + \sum_{i=1}^{n} \log[g(t_i;\boldsymbol{\xi})] - \lambda\sum_{i=1}^{n}[G(t_i;\boldsymbol{\xi})] + (\theta-1)\sum_{i=1}^{n}\log(e^{-\lambda G(t_i;\boldsymbol{\xi})} - e^{-\lambda})$$

$$- (\theta+1)\sum_{i=1}^{n} \log(1 - \alpha e^{-\lambda} - \overline{\alpha} e^{-\lambda G(t_i;\boldsymbol{\xi})}).$$

This log-likelihood function can't be solved analytically because of its complex form, but it can be maximized numerically by employing global optimization methods available with the software's R.

By taking the partial derivatives of the log-likelihood function with respect to $\theta, \alpha, \lambda$ and $\xi$ we obtain the components of the score vector $U_\rho = (U_\theta, U_\alpha, U_\lambda, U_\xi)$.

The asymptotic variance-covariance matrix of the MLEs of parameters can obtained by inverting the Fisher information matrix $I(\rho)$ which can be derived using the second partial derivatives of the log-likelihood function with respect to each parameter. The $ij^{th}$ elements of $I_n(\rho)$ are given by

$$I_{ij} = -E[\partial^2 l(\rho)/\partial \rho_i \partial \rho_j], \quad i, j = 1, 2, \cdots, 3+q.$$

The exact evaluation of the above expectations may be cumbersome. In practice one can estimate $I_n(\rho)$ by the observed Fisher's information matrix $\hat{I}_n(\hat{\rho}) = (\hat{I}_{ij})$ defined as

$$\hat{I}_{ij} \approx \left(-\partial^2 l(\rho)/\partial \rho_i \partial \rho_j\right)_{\eta=\hat{\eta}}, \quad i, j = 1, 2, \cdots, 3+q.$$

Using the general theory of MLEs under some regularity conditions on the parameters as $n \to \infty$ the asymptotic distribution of $\sqrt{n}(\hat{\rho} - \rho)$ is $N_k(0, V_n)$ where $V_n = (v_{jj}) = I_n^{-1}(\rho)$. The asymptotic behaviour remains valid if $V_n$ is replaced by $\hat{V}_n = \hat{I}^{-1}(\hat{\rho})$. Using this result large sample standard errors of $j^{th}$ parameter $\rho_j$ is given by $\sqrt{\hat{v}_{jj}}$.

**4.2 Simulation**

In this section Monte Carlo simulation study is conducted to compare the performance of the different estimators of the unknown parameters for the $GMOP-E(\theta, \alpha, \lambda, \beta)$ distribution using R program. We generate $N = 3000$ samples of size $n = 5$ to $80$ from GMOP-E distribution with true parameters values $\theta = 2, \alpha = 8, \lambda = 5, \beta = 0.5$, and calculate the bias and mean square error (MSE) of the MLEs empirically by

$$Bias_h = \frac{1}{N}\sum_{i=1}^{N}(\hat{h}_i - h) \quad \text{and} \quad MSE_h = \frac{1}{N}\sum_{i=1}^{N}(\hat{h}_i - h)^2 \quad \text{respectively (for } h = \theta, \alpha, \lambda, \beta).$$

Results of this simulation study are presented graphically in Figures 4 and 5 tells us that as the sample sizes increases the biases and MSE's approach to 0 in all cases which is consistent with the theoretical properties of the MLEs. This fact supports that the asymptotic normal distribution provides an adequate approximation to the finite sample distribution of the MLEs. The simulation study shows that the maximum likelihood method is appropriate for estimating the $GMOP-E(\theta, \alpha, \lambda, \beta)$ distribution parameters.

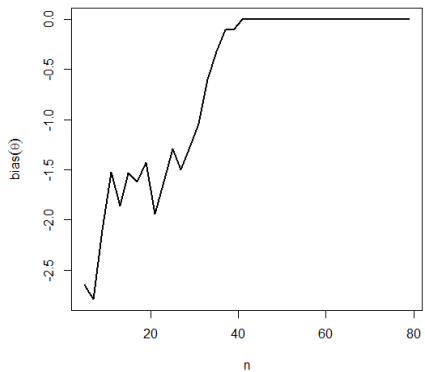
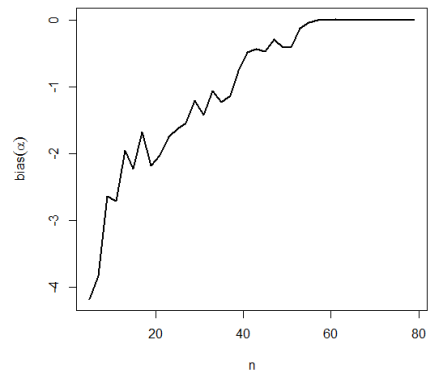

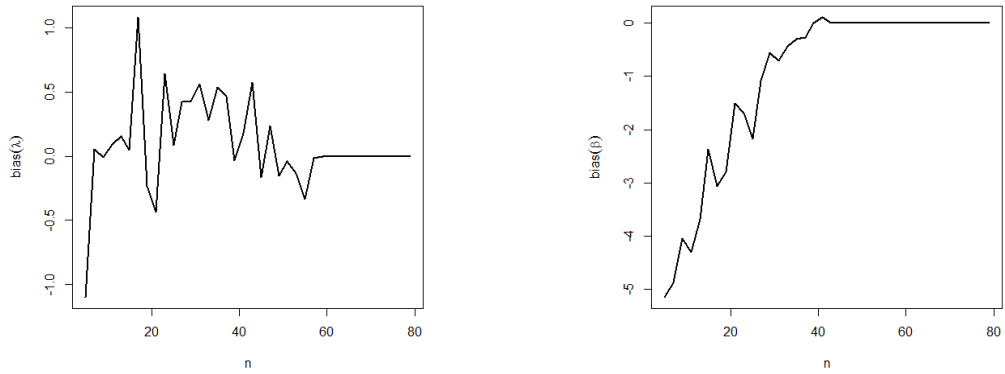

**Fig 4:** The Biases for the parameter values $\theta = 2, \alpha = 8, \lambda = 5, \beta = 0.5$ for $GMOP-E(\theta, \alpha, \lambda, \beta)$ distribution.

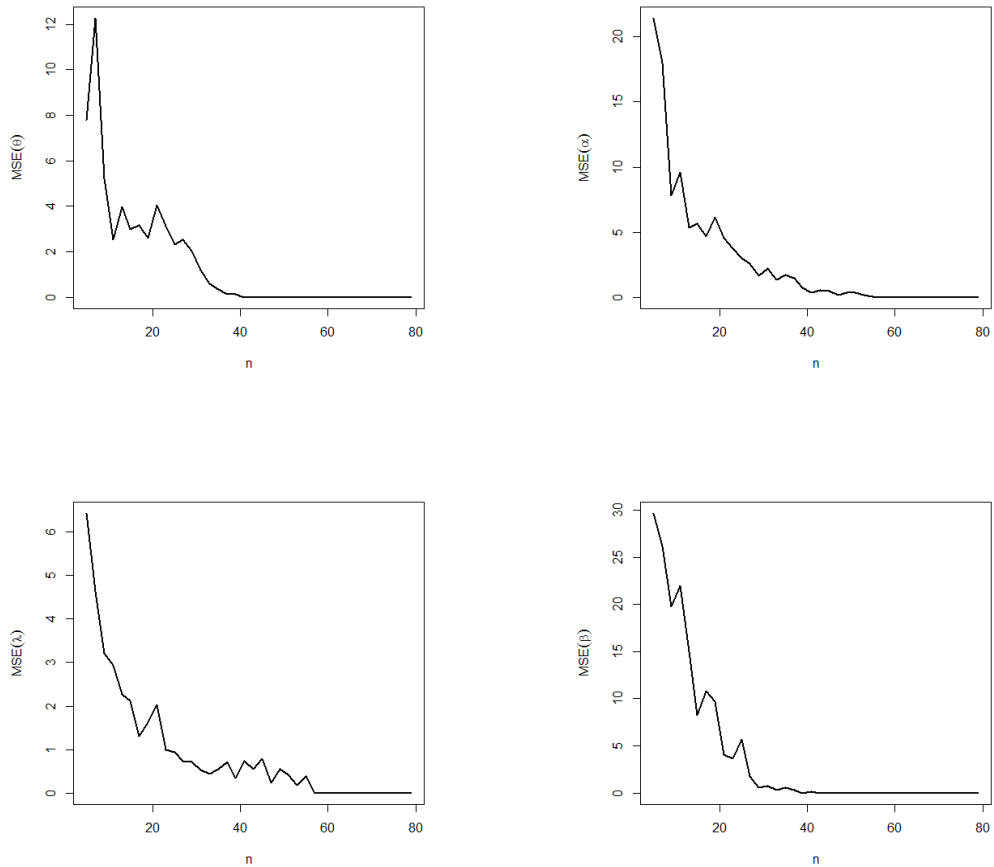

**Fig 5:** The MSEs for the parameter values $\theta = 2, \alpha = 8, \lambda = 5, \beta = 0.5$ for $GMOP-E(\theta, \alpha, \lambda, \beta)$ distribution.

## 5. Real life application

Here we consider modelling of the one failure time data set to illustrate the suitability of the $GMOP-E(\theta, \alpha, \lambda, \beta)$ distribution in comparison to some existing distributions by estimating the parameters by numerical maximization of log-likelihood functions. The data set comprises survival time of 72 guinea pigs infected with virulent tubercle bacilli, observed and reported by Bjerkedal 1960. The descriptive statistics about the data set shown in Table 3 reveal that the data set are positively skewed as expected from the nature of life time data and has higher kurtosis.

**Table 3:** Descriptive Statistics for the data set I

| Data Sets | $n$ | Min. | Mean | Median | s.d. | Skewness | Kurtosis | 1$^{st}$ Qu. | 3$^{rd}$ Qu. | Max. |
|---|---|---|---|---|---|---|---|---|---|---|
| I | 72 | 0.100 | 1.851 | 1.560 | 1.200 | 1.788 | 4.157 | 1.080 | 2.303 | 7.000 |

We have compared the $GMOP-E(\theta,\alpha,\lambda,\beta)$ distribution with some of its sub models namely, exponential (Exp), moment exponential (ME), Poisson exponential (P-E), Marshall-Olkin exponential (MO-E) (Marshall and Olkin, 1997), generalized Marshall-Olkin exponential (GMO-E) (Jayakumar and Mathew, 2008) and Marshall-Olkin Poisson exponential (MOP-E) models and also with other recently introduced models Kumaraswamy exponential (Kw-E) (Cordeiro and de Castro, 2011), Beta exponential (BE) (Eugene *et al.*, 2002), Marshall-Olkin Kumaraswamy exponential (MOKw-E) (Handique *et al.*, 2017) and Kumaraswamy Marshall-Olkin exponential (KwMO-E) (Alizadeh *et al.*, 2015), beta Poisson exponential (BP-E) (Handique *et al.*, 2020) and Kumaraswamy Poisson exponential (KwP-E) (Chakraborty *et al.*, 2020) distributions for failure time data set.

Upon fitting the best model is chosen as the one having lowest AIC (Akaike Information Criterion), BIC (Bayesian Information Criterion), CAIC (Consistent Akaike Information Criterion), and HQIC (Hannan-Quinn Information Criterion) and also, we apply formal goodness-of-fit tests in order to verify which distribution fits better to this data by considering Anderson-Darling (A), Cram′er-von Mises (W) and Kolmogorov-Smirnov (K-S) statistics to compare the fitted models. The model with minimum values for these statistics and highest *p*-value of K-S statistics could be chosen as the best model to fit the data. We have also provided the asymptotic standard errors and confidence intervals of the mles of the parameters for each competing model. For visual comparison of the best fitted density and the fitted cdf are plotted with the corresponding observed histograms and ogives in Fig. 7. These plots indicate that the proposed distributions provide a good fit to this data set.

To check the shape of the observed hazard function the total time on test (TTT) plot (see Aarset, 1987) is used. A straight diagonal line indicates constant hazard for the data set, where as a convex (concave) shape implies decreasing (increasing) hazard. The TTT plots for the data set Fig. 6 indicate that the data set has increasing hazard rate. We also provide the box plot of the data to summarise the minimum, first quartile, median, third quartile, and maximum where a box is shown from the first quartile to the third quartile with a a vertical line going through the box at the median.

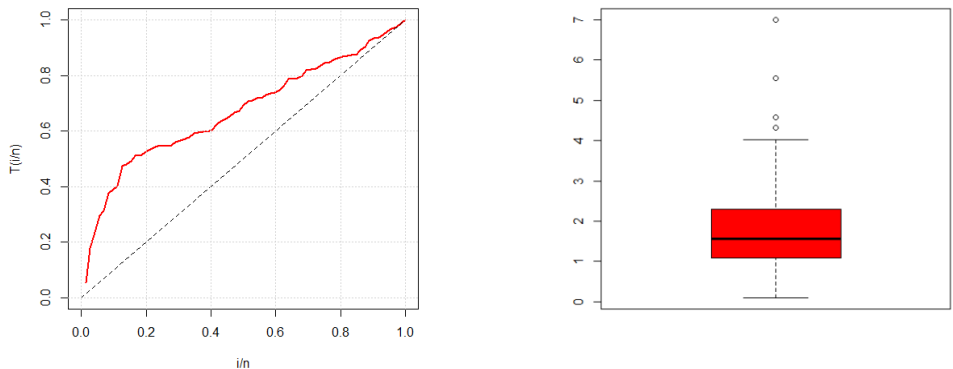

**Fig: 6** TTT and Box plot for the failure time data set

**Table 4:** MLEs, standard errors (in parentheses) values for the guinea pigs survival time's data set

| Models | $\hat{\theta}$ | $\hat{\alpha}$ | $\hat{a}$ | $\hat{b}$ | $\hat{\lambda}$ | $\hat{\beta}$ |
|---|---|---|---|---|---|---|
| Exp $(\beta)$ | --- | --- | --- | --- | --- | 0.540 (0.063) |
| ME $(\beta)$ | --- | --- | --- | --- | --- | 0.925 (0.077) |
| P-E $(\lambda, \beta)$ | --- | --- | --- | --- | -5.268 (1.069) | 1.236 (0.129) |
| MO-E $(\alpha, \beta)$ | --- | 8.778 (3.555) | --- | --- | --- | 1.379 (0.193) |
| GMO-E $(\theta, \alpha, \beta)$ | 0.179 (0.070) | 47.635 (44.901) | --- | --- | --- | 4.465 (1.327) |
| MOP-E $(\alpha, \lambda, \beta)$ | --- | 2.036 (1.071) | --- | --- | -3.366 (0.142) | 1.312 (0.093) |
| Kw-E $(a, b, \beta)$ | --- | --- | 3.304 (1.106) | 1.100 (0.764) | --- | 1.037 (0.614) |
| B-E $(a, b, \beta)$ | --- | --- | 0.807 (0.696) | 3.461 (1.003) | --- | 1.331 (0.855) |
| MOKw-E $(\alpha, a, b, \beta)$ | --- | 0.008 (0.002) | 2.716 (1.316) | 1.986 (0.784) | --- | 0.099 (0.048) |
| KwMO-E $(\alpha, a, b, \beta)$ | --- | 0.373 (0.136) | 3.478 (0.861) | 3.306 (0.779) | --- | 0.299 (1.112) |
| BP-E $(a, b, \lambda, \beta)$ | --- | --- | 3.595 (1.031) | 0.724 (1.590) | 0.014 (0.010) | 1.482 (0.516) |
| KwP-E $(a, b, \lambda, \beta)$ | --- | --- | 3.265 (0.991) | 2.658 (1.984) | 4.001 (5.670) | 0.177 (0.226) |
| GMOP-E $(\theta, \alpha, \lambda, \beta)$ | 0.333 (0.151) | 12.584 (7.696) | --- | --- | 0.054 (1.376) | 2.858 (0.959) |

**Table 5:** Log-likelihood, AIC, BIC, CAIC, HQIC, A, W and KS (*p*-value) values for the guinea pigs survival times data set

| Models | AIC | BIC | CAIC | HQIC | A | W | KS (*p*-value) |
|---|---|---|---|---|---|---|---|
| Exp ($\beta$) | 234.63 | 236.91 | 234.68 | 235.54 | 6.53 | 1.25 | 0.27 (0.06) |
| ME ($\beta$) | 210.40 | 212.68 | 210.45 | 211.30 | 1.52 | 0.25 | 0.14 (0.13) |
| P-E ($\lambda, \beta$) | 209.86 | 214.42 | 210.03 | 211.66 | 0.98 | 0.19 | 0.09 (0.66) |
| MO-E ($\alpha, \beta$) | 210.36 | 214.92 | 210.53 | 212.16 | 1.18 | 0.17 | 0.10 (0.43) |
| GMO-E ($\theta, \alpha, \beta$) | 210.54 | 217.38 | 210.89 | 213.24 | 1.02 | 0.16 | 0.09 (0.51) |
| MOP-E ($\alpha, \lambda, \beta$) | 208.32 | 215.15 | 208.67 | 211.04 | 0.96 | 0.17 | 0.09 (0.56) |
| Kw-E ($a, b, \beta$) | 209.42 | 216.24 | 209.77 | 212.12 | 0.74 | 0.11 | 0.08 (0.50) |
| B-E ($a, b, \beta$) | 207.38 | 214.22 | 207.73 | 210.08 | 0.98 | 0.15 | 0.11 (0.34) |
| MOKw-E ($\alpha, a, b, \beta$) | 209.44 | 218.56 | 210.04 | 213.04 | 0.79 | 0.12 | 0.10 (0.44) |
| KwMO-E ($\alpha, a, b, \beta$) | 207.82 | 216.94 | 208.42 | 211.42 | 0.61 | 0.11 | 0.08 (0.73) |
| BP-E ($a, b, \lambda, \beta$) | 205.42 | 214.50 | 206.02 | 209.02 | 0.55 | 0.08 | 0.09 (0.81) |
| KwP-E ($a, b, \lambda, \beta$) | 206.63 | 215.74 | 207.23 | 210.26 | 0.48 | 0.07 | 0.09 (0.79) |
| **GMOP-E** ($\theta, \alpha, \lambda, \beta$) | **204.24** | **213.36** | **204.83** | **207.84** | **0.44** | **0.04** | **0.07 (0.83)** |

In the Tables 4 and 5 the MLEs with standard errors of the parameters for all the fitted models along with their AIC, BIC, CAIC, HQIC, A, W and K-S statistic with *p*-value for the six sub models for the failure time data set are presented respectively. From these findings, it is evident that the GMOP-E distribution with lowest value of AIC, BIC, CAIC, HQIC, A, W and highest *p*-value of K-S statistics is not only a better model than all the sub models Exp, ME, P-E, MO-E, GMO-E, MOP-E but also the better than the most of the recently introduced three or four parameters models namely Kw-E, B-E, MOKw-E, KwMO-E, BP-E and KwP-E. These findings are further validated from the plots of fitted density with histogram of the observed data and fitted cdf with ogive of observed data in Figure 7. These plots clearly indicate that the proposed distribution provide closest fit to the data set considered here.

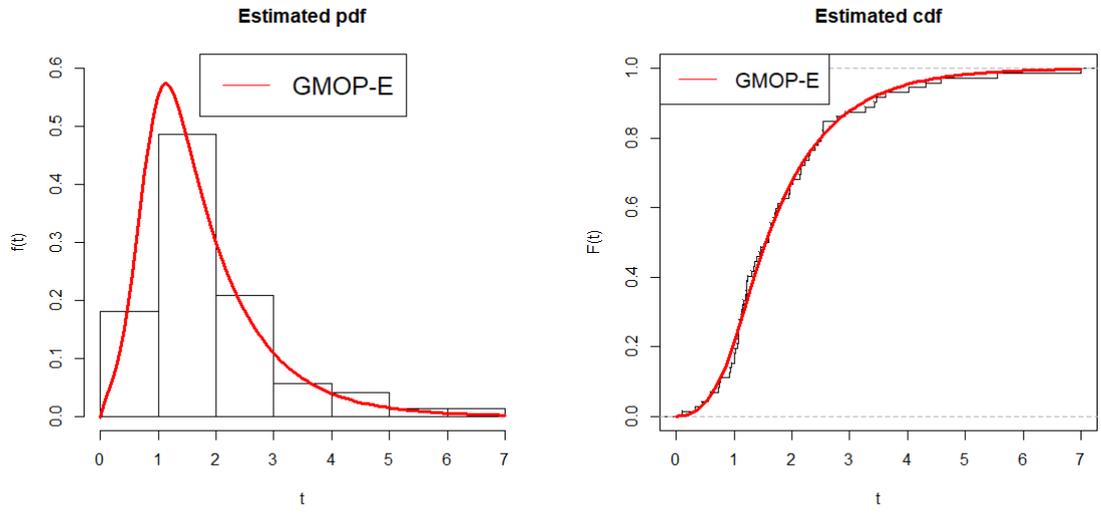

**Fig: 7** Plots of the observed histogram and estimated pdf on left and on right the observed ogive and estimated cdf for failure time data set for the GMOP-E model

## 6. Conclusions

In this work, we propose a new family of continuous distributions called the *Generalized Marshall-Olkin Poisson -G* family of distributions. Several new models can be generated by considering special distributions for G. We demonstrate that the pdf of any GMOPG distribution can be expressed as a linear combination of exponentiated-G density functions, which result allowed us to derive some of its mathematical and statistical properties such as moment generating function, order statistics, probability weighted moments and Rényi entropy. The estimations of the model parameters are obtained by maximum likelihood method. One application of the proposed family empirically prove its flexibility to model the real data sets, in particular we verified that a special case of the GMOPG family can provide better fits than its sub models and other models generated from well-known families.